\documentclass{elsarticle}

\newcommand{\be}{\begin{eqnarray}}
\newcommand{\ee}{\end{eqnarray}}

\newcommand{\short}{\mbox{Short}}
\newcommand{\cut}{\mbox{Cut}}
\newcommand{\m}{\setminus}

\newtheorem{theorem}{Theorem}

\newtheorem{lemma}{Lemma}

\usepackage{graphicx}
\usepackage{amssymb}

\begin{document}

\title{Threat, support and dead edges in the Shannon game}

\author{Andrew M. Steane}
\ead{a.steane@physics.ox.ac.uk}
\address{
Department of Atomic and Laser Physics, Clarendon Laboratory,\\
Parks Road, Oxford, OX1 3PU, England.
}

\date{\today}

\begin{abstract}
The notions of {\em captured/lost vertices} and {\em dead edges} in the Shannon 
game (Shannon switching game on nodes) are examined using graph theory. 
Simple methods are presented for identifying some dead edges and
some captured sets of vertices, thus simplifying the (computationally hard)
problem of analyzing the game.
\end{abstract}


\maketitle

\section{Introduction}

In the {\em Shannon game} two vertices in a simple graph are designated as terminals,
and one player aims to link the terminals, the other to prevent
this, by colouring vertices. 
This is of interest both in its own right as a mathematical
problem in graph theory and because of its connection to the theory of
robust networks and fault-tolerance. 
An understanding of the game can
elucidate more general problems, such as the Shannon game with
more terminals (discussed below) or with terminals that can
change. Various hard computational problems arise in aspects of
the game, see \cite{06Hayward} for a review.

It was shown by Bj\"ornsson {\em et al.} \cite{06Bj} that analysis of the Shannon game can be simplified by
considering local games (played on a subgraph) as long as the complete neighbourhood of the
subgraph is considered, yielding a type of Shannon game with multiple terminals, called
a multi-Shannon game (defined in \cite{06Bj} and below). They introduced a concept of
{\em domination}, relating to the winner of a multi-Shannon game, and presented a
set of local patterns in the game {\em Hex} where empty cells could be filled
in without changing the outcome of the game.

In this paper we introduce concepts related to domination (in the sense of Bj\"ornsson
{\em et al.}), but which are defined in terms of easily identifiable properties of a graph.
This simplifies the analysis and permits a number of theorems concerning multi-Shannon
games to be easily stated and proved. All the local {\em Hex} patterns
considered in \cite{06Bj} can be detected through simple graph properties such as
counting small walks and triangles. We also discuss the concept of a dead edge, defined
in section \ref{s.dead}, and present theorems to help identify dead edges.
Some further theorems allowing won or lost
patterns to be detected are also obtained. 

\section{Game definition}

The Shannon game is played on any simple
graph (i.e. no self-loops and no multiple edges),
with two of the vertices designated terminals. One player (called \short)
seeks to join the terminals, the other  (called \cut) seeks to prevent
this. Shannon originally proposed colouring edges of the graph, i.e. in each
turn \short~can colour one edge, and \cut~can delete one uncoloured edge. This
version is known as the {\em Shannon switching game}; it has a polynomial-time
algorithmic solution \cite{64Lehman}. In this paper we consider only
the Shannon switching game on nodes, known as the {\em Shannon game}. In the
Shannon game, in each turn \short~can colour one
vertex, and \cut~can delete one uncoloured vertex. \short~has won if there
exists a path between the terminals consisting only of coloured vertices; \cut~
has won if the graph is disconnected with the terminals in different components.
It is obvious that a draw is not possible.

Instead of deleting vertices \cut~may
equivalently colour vertices a different colour to \short~(e.g. \short~is
black and \cut~is white), then \cut~has won if there is a
cutset separating the terminals consisting of all white vertices. We will
call this the {\em graph colouring model} of the game.

The game {\em Hex} can be considered as a (beautiful and fascinating) example of
the Shannon game, played on a finite planar graph corresponding to an
$m \times m$ honey-comb
lattice of hexagonal cells (vertices in the graph are adjacent
when their corresponding cells are
adjacent), with two further vertices (the terminals) adjacent to all vertices
on each of two opposite sides of the lattice \cite{06Hayward}. If we choose to define the
rules of Hex such that
one player has to connect the designated sides and the other has to prevent this, then
the correspondence to the Shannon game is direct and obvious. However, the rules
of Hex are conventionally stated another way: both players colour vertices
(with different colours) and one player has to connect one pair of sides, the other
player the other pair of sides. That the latter version (standard Hex rules)
is equivalent to the former (Shannon game)---i.e. that joining one pair of sides
in Hex is equivalent to separating the other pair of sides---was shown by Beck \cite{69Beck} and
Gale \cite{79Gale}, see also \cite{06Hayward}.

The Shannon game as described above, where in each turn \short~colours a
vertex, can also be modelled another way, in which whenever \short~picks
a vertex $v$, edges are added to the graph between
all pairs of neighbours of $v$ which are not already adjacent, and
then $v$ is deleted.
We refer to this as the
{\em reduced graph model}. It is easy to see that the
graph colouring and reduced graph models are equivalent. Most
of the proofs in this paper use the reduced graph model.

The reduced graph model shows that any position in the course of a
Shannon game is equivalent to an opening position on a smaller graph
(i.e. one with fewer vertices; the number of edges may grow or shrink). Such
a position is uniquely specified by the reduced graph, the terminal pair
designation and identity of the player who is to move next:
$(G,\{t_1,t_2\},P)$. In this paper the Roman font symbol $G$ will always refer
to a reduced graph, and we will usually suppress the terminal designation
in order to reduce clutter (the terminal designation does not change during
the game). For convenience we assign binary values to the players,
$P=1$ for \short~and $P=0$ for \cut, so if a given player is $P$, then the other
player is $\neg P$.

\begin{lemma} \label{lm:1}
For every position in the Shannon game, one of the
players has a winning strategy.
\end{lemma}

This is well known and follows from elementary game theory, but it is essential
in the following so we exhibit a proof, by induction: Let $(G_n,P)$ describe
the position when the reduced
graph has $n$ vertices, the player to move is $P$, and we suppressed the terminal
designation. We
assume that for every possible position of size $n-1$ which is reachable from
$(G_n,P)$, one of the players (not necessarily the same for all the graphs) has a
winning strategy. If among this set there is a position $(G_{n-1},\neg P)$
for which $P$ has a winning strategy, then $P$ has a winning strategy
for $(G_n,P)$. If not, then all the $\{(G_{n-1},\neg P)\}$
reachable from $(G_n,P)$ give a winning
strategy for $\neg P$, from which it follows that $\neg P$ has a winning
strategy for $(G_n,P)$. This establishes the inductive step.
It remains to prove the case for some low $n$,
and this follows immediately from the fact that $n$ falls as the game
proceeds, and ultimately the game is won or lost (not drawn).

The analysis of the Shannon game has two main aims: to identify the winner
for a general position, and to find a winning strategy.

Determining the winner of a general Hex position has been shown to be
PSPACE-complete \cite{81Reisch}, and therefore so is the solution of
the Shannon game\footnote{This does not rule out that to determine the winner
of some specific class of positions, such as the opening position in
$m \times m$ Hex, may be easier.}.

\subsection{Notation and terminology}

We will use the term `graph' to signify a simple graph with two special vertices
(the terminals), i.e. $G = (V,E,\{t_1,t_2\})$ is a simple graph
with vertex set $V$, edge set $E$, terminal set $\{t_1,t_2\} \subset V$.
A {\em position} $(G,P)$ is
a graph (having designated terminals) with an indication $P$ of which player
is to play next.

A {\em clique} is a set of vertices each adjacent to all the others
in the set. In standard graph terminology the term means a maximal
such set. In this paper it will not be important whether the set
is maximal or non-maximal (i.e. a subset of a maximal clique). We
will use `clique' to mean a vertex set which is either a
non-maximal clique or a (maximal) clique.

$G-v$ means the graph obtained from $G$ by deleting vertex $v$. $G
\ast v$ means the graph obtained from $G$ by shorting vertex $v$,
that is, adding edges between all non-adjacent pairs in the
neighbourhood of $v$ (converting the neighbourhood into a clique),
and then deleting $v$. Note that in both cases the vertex is
deleted.  These two operations commute when applied to two different
vertices (this is obvious from the
graph colouring model). It is not possible to apply them successively to
the same vertex. For an edge $e$, $G+e$ means the graph
obtained by adding edge $e$, $G-e$ means the graph obtained by
deleting edge $e$.

If $S$ is a set vertices, then $G-S$ means the graph obtained by deleting
the members of $S$ and $G * S$ means the graph obtained by shorting the members of $S$. 

Terminals may never be deleted or shorted.

The open neighbourhood of a vertex $v$ is denoted $\Gamma(v)$. The
closed neighbourhood is denoted $\Gamma[v] \equiv \Gamma(v) \cup
\{v\}$. Define the neighbourhood $\Gamma(S)$ of a set of vertices
$S$ to be the set of vertices not in $S$ but adjacent to a member
of $S$: $\Gamma(S) \equiv ( \cup_{v \in S} \Gamma(v) ) \m S$,
where $\m$ is the set minus operation.

If $a \in S$ then we use the shorthand $S-a \equiv S \m \{a\}$.

The triangle number $T(v)$ is the number of triangles containing
vertex $v$.

The {\em win-value} $W(G,P)$ of a position is the identity of the player
holding a winning strategy.

\section{Multi-Shannon game: domination, capture and loss}  \label{s:multi}

The {\em multi-Shannon game} is introduced in \cite{06Bj} in order to define
concepts and methods useful to solving the Shannon game.
In this section we present the game definition,  
quote the main result of \cite{06Bj}, and present some
new results.

If $S$ is a set of vertices in a graph $G$ then a multi-Shannon game is played on
the induced subgraph $G(S \cup \Gamma(S))$ with the vertices $S$ as `playing area' and
the vertices $\Gamma(S)$ designated as `terminals'. There may be any number of
terminals.

To describe the game it is convenient to use both the graph colouring model and the
reduced graph. \short\ seeks to play in such a way that at the end of the game, when
all vertices in $S$ have been coloured, the reduced graph (which, at the end, consists only of
terminals and possibly edges between them) is the same as would be obtained if all 
of $S$ had been shorted. \cut\ seeks to play in such a way that the outcome 
(i.e. the final reduced graph among terminals) is the same as if all of $S$ had
been cut. If neither eventuality occurs, then the game is drawn. (If both occur, because,
for example, the terminals were all mutually adjacent at the start, then both players
may declare a win; the vertices in $S$ are then `dead', a concept we will discuss later).

The above rule definition clears up a slight ambiguity in \cite{06Bj} for the case where
the graph was not connected at the outset.
If the graph is connected then a \short\ win turns the reduced graph terminals into a clique.
If the terminals initially had no edges between them
then a \cut\ win turns the reduced graph terminals into a set of isolated vertices.
In the graph colouring model, \short\ is trying to make \short-coloured paths 
between all connected but non-adjacent terminal
pairs, and \cut\ is trying to occupy cutsets that separate all non-adjacent terminal pairs.

The multi-Shannon game is useful because it offers a way to simplify the
analysis of the regular Shannon game, by reducing part of the complete
game to a local `battle', as follows. Consider any set $S$ of non-terminal
vertices of the reduced graph in a Shannon game. Let
${\cal M}_{G,S}$ be the multi-Shannon game played on the
subgraph of $G$ induced by $S \cup \Gamma(S)$ with terminals $\Gamma(S)$.

{\bf Definition} $S$ is {\em captured} in $G$ if \short~has a
second-player winning strategy for ${\cal M}_{G,S}$. $S$ is {\em
lost} in $G$ if \cut~has a second-player winning strategy for
${\cal M}_{G,S}$.

\begin{theorem} \label{BHJB}
(Bj\"ornsson, Hayward, Johanson, van Rijswijck)
If $S$ is captured, then $W(G,x) = W(G \ast S,x) \forall x$.
If $S$ is lost, then $W(G,x) = W(G - S,x) \forall x$.
\end{theorem}

Proof: see \cite{06Bj}.
Our terminology `captured' and `lost' reflects the
point of view of \short. The same concept is referred to in \cite{06Bj} as
`captured by short' and `captured by cut'.

\begin{figure}[ht!]
\centerline{\resizebox{!}{0.1\textwidth}{\includegraphics{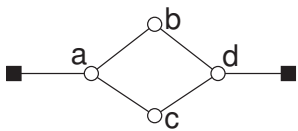}}}
\caption{Examples of capture and loss. The filled squares are terminals. 
$\{b,c\}$ is captured but $\{a,b,c,d\}$ is lost.}
\label{f:captureloss}
\end{figure}

Observe that in the first case, to solve $G$, it suffices to solve
$G \ast S$, and in the second case, to solve $G$, it suffices to
solve $G - S$. In either case the resulting problem is usually
simpler because the reduced graph has fewer vertices. Therefore
when vertices are captured (lost) one may as well go ahead and
short (cut) them as a `free move' in order to clarify matters.
This is called `filling in' in discussions of Hex \cite{06Hayward}. To be
clear in the following one should note that to be captured or lost
is strictly a property of a set of vertices (in a given graph) rather than of the
members of the set. For example a captured set can be a subset of a
lost set, and vice versa, see figure \ref{f:captureloss} for
examples.

{\bf Definition} $S$ is {\em $P$-dominated} if $P$ has a first-player
winning strategy for ${\cal M}_{G,S}$. For any initial move $v$ in such a
strategy, we say that $v$ {\em $P$-dominates} $S$.

\begin{figure}[ht!]
\centerline{\resizebox{!}{0.125\textwidth}{\includegraphics{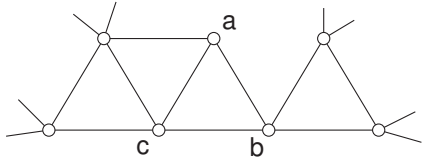}}}
\caption{Examples of cut-domination.
The diagram shows part of a larger graph; edges attached to only one vertex are used
to signify that the vertex may have further neighbours which are irrelevant to the
point under discussion. In this example,
$b$ cut-dominates $\{a,b\}$ because after deleting $b$ the neighbourhood of $a$ becomes
a clique. $a$ cut-dominates $\{a,b,c\}$ because if $a$ is deleted then $b$ and $c$
become a lost pair (c.f. figure \protect\ref{fig:pair_example}(a)).
However, $a$ does not cut-dominate $\{a,b\}$ or $\{a,c\}$.
(The edge $ac$ is dead, see theorem \protect\ref{th:vwdead}).
} \label{fig:dominate}
\end{figure}

The terminology here is precisely as in \cite{06Bj}. Domination
is a property of a set as a whole:
if a set $S$ is $P$-dominated it does not necessarily follow
that a subset $U \subset S$ is $P$-dominated. For example figure \ref{fig:dominate} shows
a case where $a$ cut-dominates $\{a,b,c\}$ but $a$ does not cut-dominate $\{a,b\}$
or $\{a,c\}$.

{\bf Simple examples} .
If $G-a \ast b-c = G-a-b\ast c = G-a-b-c$ then $a$ cut-dominates $\{a,b,c\}$
(c.f. figure \ref{fig:dominate});
if  $G\ast a \ast b-c = G\ast a - b\ast c = G\ast a \ast b \ast c$
then $a$ short-dominates $\{a,b,c\}$.

\begin{theorem} \label{th:doubledom}
If both vertices in a two-vertex set $P$-dominate the set, then
both are captured, lost respectively for $P$=\short, \cut.
\end{theorem}

(This property was noted in \cite{06Bj}).
Proof: in the multi-Shannon game played on these two vertices,
after the opening move by $\neg P$, there is only
one response for $P$ and then the game is over, with a final position exactly
the same as if $P$ had opened (with a winning move) and $\neg P$ replied. $\Box$

\begin{figure}[ht!]
\centerline{\resizebox{!}{0.15\textwidth}{\includegraphics{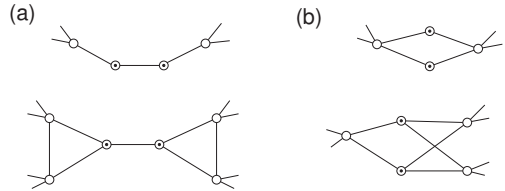}}}
\caption{Examples of mutually dominating pairs.
(a) The dotted vertices cut-dominate each other and therefore are lost.
(b) The dotted vertices short-dominate each other and therefore are captured.
} \label{fig:pair_example}
\end{figure}

{\bf Simple examples} (See figure \ref{fig:pair_example}).
Two neighbouring vertices both of degree 2 are lost, unless
one or both is terminal. For two
vertices of degree 2 having the same neighbours,
if neither is terminal then both are captured; if one is terminal then the
other is dead (see later); if both are terminal then the game is won by \short.

N.B. This theorem is concerned with pairs not larger sets. It states for example that if
$a$ cut-dominates $\{a,b\}$ and $b$ cut-dominates $\{a,b\}$ then $a$ and $b$ are
lost. However, if $a$ cut-dominates $\{a,b,c\}$ and $b$ cut-dominates $\{a,b\}$ it
does not necessarily follow that $\{a,b\}$ is lost: figure \ref{fig:dominate} gives a
counter-example.

\begin{theorem} \label{th:threedom}
If each of two vertices in a three-vertex set $P$-dominate the set, then so does the third.
\end{theorem}

Proof: let $1$ and $0$ stand for moves by $P$ and $\neg P$.
There are three possible outcomes after three moves in a three-vertex
multi-Shannon game when $P$ goes first: $110,\, 101,\, 011$ (in an obvious notation).
If the first vertex $P$-dominates the set then both
$110$ and $101$ must be $P$-wins because $\neg P$ can force either of these
positions after
a $P$ opening at the first vertex. Similarly, if the
second vertex $P$-dominates the set then both
$110$ and $011$ must be $P$-wins. Hence $101$ and $011$ are both
$P$-wins, and therefore the third vertex also $P$-dominates the set.
$\Box$

For example, in figure \ref{fig:dominate} $a$ and $b$ both cut-dominate $\{a,b,c\}$,
therefore so does $c$.

\begin{theorem} \label{th:fourdom}
If each of two vertices in a four-vertex set $P$-dominate the set, but occupying
them both is not a winning strategy for $P$, then $P$ 
has a second-player winning strategy for the multi-Shannon game on this set.
\end{theorem}

\begin{figure}[ht!]
\centerline{\resizebox{!}{0.15\textwidth}{\includegraphics{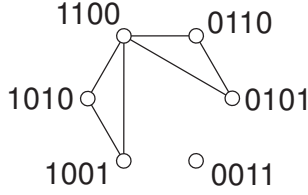}}}
\caption{A graph indicating possible outcomes of a game on a 4-vertex set considered
in theorem \protect\ref{th:fourdom}.} 
\label{f.ms4pos}
\end{figure}

Proof: let the vertices be labeled $a,b,c,d$ with $a$ and $b$ $P$-dominating the set.
There are six possible final positions, namely $1100, 1010, 1001, 0110, 0101, 0011$,
where we list the vertex colours in the order $abcd$. Suppose first that $P$
moves first, and occupies vertex $a$. Since $a$ $P$-dominates the set, this is a winning move,
which implies that for each of the three replies available to $\neg P$, $P$ can
still obtain a won final position. This means that $P$ must win at least one
of $1010,1001$ (the available outcomes if $\neg P$ replies at $b$), and at least one of
$1100,1001$ and at least one of $1100,1010$. These pairs are indicated by three of the
edges of the graph shown in figure \ref{f.ms4pos}.  Now suppose that $P$ moves
first and occupies vertex $b$. Since $b$ $P$-dominates the set, this is also a winning move
and by similar reasoning we conclude that $P$ wins at least one of each of the
three pairings indicated by the remaining edges in figure \ref{f.ms4pos}.

Let us call the graph shown in figure \ref{f.ms4pos} the {\em win-graph}. 
Next consider the win-values of each of the final positions. 
For each edge in the win-graph, at least one vertex to which it is incident is a won
position for $P$. In other words, the won positions form a vertex-cover of the win-graph.
However, the conditions of the theorem assert that $1100$ is not winning for $P$. It follows
that $1010,1001,0110,0101$ are all won by $P$. Therefore the only way for $\neg P$ to 
win or draw is to realise one of $1100$ or $0011$. However, after any opening move by
$\neg P$, the reply by $P$ can guarantee to avoid both of those outcomes. Hence $P$
must win. $\Box$

The general relationship between capture/loss and domination is obvious:
$P$ has a second-player winning strategy if and only if $P$ has a first-player winning strategy
after any opening move by $\neg P$. For the case of \short\ this general idea can be simplified,
as follows.

\begin{theorem} \label{th:winmulti}
$S$ is captured if and only if
\short~has 1st-player winning strategies for ${\cal M}_{G,S}$
on each of a set of subsets of $S$ having empty total intersection.
\end{theorem}

By ``subsets having empty
total intersection" we mean $\{U_i\} : U_i \subset S, U_1 \cap U_2 \cap \ldots U_m = \emptyset$,
where $m$ is the number of subsets. A winning strategy for $P$ ``on a subset" means 
that the strategy involves play (by either player) only within the subset.
To prove such a strategy on a subset for \short, it suffices
to consider the subgraph induced by the subset plus $\Gamma(S)$ (i.e. the 
terminals of 
${\cal M}_{G,S}$, not the neighbourhood of the subset under consideraton), and
this means the theorem serves as a
useful way to break down the task of finding 2nd-player winning strategies for \short.
The subsets useful for \short\ are those which connect to all the terminals and which 
are well-connected internally. We will explain in a moment why there is no
straightforwardly analogous theorem for \cut.

Proof. ($\Leftarrow$)
Cutting vertices outside a given subset does nothing to the graph induced by that 
subset, and the empty total intersection implies that after the opening move by \cut, 
at least one \short-winning subset has not been intruded upon. \short~plays
the 1st-player winning strategy for any one such subset.
($\Rightarrow$) If \short~has a 2nd-player winning strategy, then by definition he has
a 1st-player winning strategy after any opening move $v$ by \cut, and that move did not 
affect the subgraph induced by $U_v = S-v$. Therefore the set of subsets
$\{U_v\} \forall v \in S$ is one possible set satisfying the theorem. $\Box$

When one tries to break down the task of \cut~into sub-games, the situation
is less straightforward, because ruling out the existence of any path is more 
complicated than establishing the existence of one path. We already noted that \cut\ has a 2nd-
player winning strategy if and only if he has a 1st-player winning 
strategy after any opening move by \short, but one would like to identify a way
to reduce the problem further. One possibility is to look for a smaller
multi-Shannon game (a sub-game of the one under consideration)
which has a 1st-player strategy for \cut\ to win outright, but
there does not necessarily exist any such sub-game. One cannot simply adapt
theorem \ref{th:winmulti} to \cut's case,
for two main reasons.
First, to prove
that \cut~can prevent paths useful to \short~from forming within a
given vertex set $U$, it is not
sufficient to consider the induced subgraph $G(U \cup \Gamma(S))$,
because shorting vertices outside of $U$ can add edges to $G(U)$.
Secondly, if it is to be
sufficient for \cut~to have first-player wins on each
of a set of subsets, then the subsets must be arranged in
series rather than in parallel, that is,
nested in such a way that every path between non-neighbour terminals passes through
sufficiently many subsets.
However, the property of whether a given subset $U \subset S$ is thus in series with the
rest is not a property of $U$ alone: it is a property of the whole set\footnote{An
exception is where a subset contains the neighbourhood of all terminals,
or of all terminals except a group of terminals forming a clique.}.
This is in contrast with the subsets sought by \short: the property
that a given subset $U$ connects all terminals $T = \Gamma(S)$
in a multi-Shannon game (i.e. $T \subset \Gamma(U)$) is a property of
the induced subgraph $G(U \cup T)$ alone, without regard to the rest of the graph.

Theorem \ref{th:winmulti} is reminiscent of the `{\sc or}-rule' recognized
by all early exponents of the game of Hex and investigated more 
fully by Anshelevich \cite{00Anshelevich,02Anshelevich}.
The standard `{\sc or}-rule' is concerned with `virtual connections' or `links' between
a given pair of vertices; its generalization to the case of
connections between any number of vertices is discussed in the appendix.

%A win ``on each of two distinct subsets" means that there are two 1st-player
%winning strategies for $P$, each requiring $P$ to play on only
%one of the subsets. This can happen for \short~when
%$S = U \cup V$ with $\Gamma(U) \m V = \Gamma(V) \m U = \Gamma(S)$,
%and for \cut~when $S=U \cup V$ with $U$ and $V$ nested such that all paths
%between non-adjacent terminals of $S$ pass through both $U$ and $V$
%(see figure \ref{fig:mgame_abcd} for examples).
%The proof is obvious: whichever subset $\neg P$~opens on, $P$~plays the
%winning strategy on the other, and therefore has a 2nd-player win overall.
%The statement of theorem \ref{th:winmulti} is not fully
%symmetric between \short~and \cut~because the graph induced by a given
%subset plus terminals can not be modified by deleting vertices
%outside the subset, but it can be modified by shorting vertices outside
%the subset: this can introduce new edges\footnote{For example, consider a
%degree 4 vertex adjacent to four degree 2 vertices, each adjacent to one of
%4 pendant terminals.}.

Although some of the theorems and proofs relating to capture and
loss are similar, nevertheless the goals of \short~and \cut~are not
symmetric when playing the Shannon game on a general graph (as opposed to Hex),
and this results in the asymmetry just discussed, and the asymmetry between
parts (i) and (ii) of theorem
\ref{th:threat_support} below. It is interesting to ask
whether a general Shannon game can be analysed as if \short~and \cut~were
both trying to connect different pairs of terminals, as is the case for Hex.
This would be the case if we can attach two further vertices to the graph
(to serve as \cut's terminals), such that all cutsets winning the Shannon
game for \cut~are paths between these two extra vertices, and all paths
winning the game for \short~are cutsets separating the two extra vertices.
Intuitively it seems that this should always be possible for planar graphs
which can be embedded in the plane with the terminals at the outside, and for graphs
obtainable from such a planar graph by shorting.
It may or may not be possible for other graphs. When it is possible,
one can further analyse the game by maintaining two `dual' graphs, in which
a short operation in one graph is a deletion in the other, and vice-versa.
However, we will not pursue this idea further here.

\section{Threat and support}

The notions of capture, loss and domination are useful because
they can reduce part of the Shannon game to smaller sub-problems.
However, this is only useful if we can recognize and solve those
smaller problems efficiently. In this and the next section we
consider this issue. We introduce concepts of {\em support} and
{\em threat}, which are like domination but which are defined in
terms of easily recognizable properties of the graph itself, not
in terms of the winner of a multi-Shannon game. This approach
allows some useful theorems to be stated and proved.

{\bf Definition} For terminal-free, non-intersecting vertex sets
$A,B$, $A$ {\em threatens} $B$ in $G$ if $G - A \ast B = G - A -
B$. $A$ {\em supports} $B$ in $G$ if $G \ast A - B = G \ast A \ast
B$. (Terminals can neither threaten or support nor be threatened
or supported).

The idea of the first definition is that if \short~shorts one or
more vertices in $B$, then cutting $A$ (for example by cutting a
vertex that cut-dominates $A$ and then deleting lost vertices)
will yield a graph that is identical to one which would exist if
$B$ had been deleted. Similarly, in the second definition, if
\cut~cuts a vertex in $B$, then shorting $A$ will `mend the
damage'. Threat and support are stronger properties than
domination: if vertex $v$ supports or threatens vertex-set $B$ 
then $v$ dominates $\{v\} \cup B$ but the converse is not 
necessarily true (a vertex which dominates a set does not necessarily
support or threaten the set).
Threat and support are less subtle than domination,
but easier to detect by examining local graph properties.

If $A$ is cut-dominated and $A$ threatens $B$, then
$A \cup B$ is cut-dominated, and furthermore a \short~move at any
vertex in $B$ is locally losing (see lemma \ref{lm:robust}(i)).
If $A$ is short-dominated and $A$
supports $B$ then $A \cup B$ is short-dominated and furthermore a
\cut~move at any vertex in $B$ is locally losing
(see lemma \ref{lm:robust}(iii)). 

\begin{theorem} \label{th:domts}
Domination of a two-vertex set is identical to threat or support. That is,
(i) $v$ \cut-dominates a two-vertex set $\{v,w\}$ in $G$ if and only if
$G - v \ast w = G-v-w$; and
(ii) $v$ \short-dominates a two-vertex set $\{v,w\}$ in $G$ if and only if
$G \ast v - w = G \ast v \ast w$.
\end{theorem}

Corollary: A pair of mutually supporting vertices are captured; a
pair of mutually threatening vertices are lost.

Proof:  Let $T = \Gamma(\{v,w\})$ before $v$ or $w$ are cut or
shorted. (i) ($\Rightarrow$) If $v$ \cut-dominates $\{v,w\}$ then
after cutting $v$, shorting $w$ must introduce no new edges among
the vertices of $T$ (if it did then the multi-Shannon game cannot
have been won by \cut). However after cutting $v$, the
neighbourhood of $w$ is a subset of $T$. Therefore shorting $w$
must not now introduce edges between any neighbours of $w$,
therefore shorting $w$ now has the same effect on the graph as
cutting $w$. ($\Leftarrow$) is obvious. $\Box$ (ii)
($\Rightarrow$) By the definition of the operations $\{\ast,-\}$,
$G\ast v - w$ can only differ from $G\ast v \ast w$ by having
fewer edges between members of $T$. However, if $v$
\short-dominates $\{v,w\}$ then $G \ast v - w$ is a \short-win for
the multi-Shannon game on $\{v,w\}$, and the definition of a
\short~win is that the adjacency among $T$ is the same as if both
$v$ and $w$ were shorted. It follows that, in such conditions,
there can be no difference between $G\ast v - w$ and $G\ast v \ast
w$. ($\Leftarrow$) is obvious. $\Box$
The corollary follows immediately by combining this theorem
with theorem \ref{th:doubledom}.

We now consider two issues: how to identify threats and supports,
and their robustness as the game is played. First we need a lemma
concerning sub-graphs with clique neighbourhoods (or cutsets). 
Parts (i,ii) of the lemma
introduce the main property, and part
(iii) shows that once a sub-graph has a clique neighbourhood,
further play preserves that property.
By a `connected component of $S$' we mean a vertex set that forms a maximal
connected component of the induced subgraph $G(S)$.

\begin{lemma} \label{lm:clique}
(i) For a set of vertices $S$ whose induced subgraph $G(S)$ is connected,
$G \ast S = G - S$ if and only if $\Gamma(S)$ is a clique. \\
(ii) For any set of vertices $S$ in $G$, $G \ast S = G - S$ if and
only if $G \ast S_i = G - S_i$
for each connected component $S_i$ of $S$. \\
(iii) If $G \ast S = G - S$ for a set of vertices $S$,
then for any vertex set $V$, $(G-V) \ast (S \m V) = (G-V) - (S \m V)$
and $(G \ast V) \ast (S \m V) = (G\ast V) - (S \m V)$.
\end{lemma}

Proof: (i) ($\Leftarrow$) If $\Gamma(S)$ is a clique then the
operations $\ast S$ and $-S$ are strictly identical.
($\Rightarrow$) The operation $\ast$ consists of edge insertion
followed by the operation $-$. The only edges affected by the
operation $-S$ are those incident on members of $S$. Therefore if
$G \ast S=G-S$ it must be that the operation $\ast S$ inserts no
edges between non-members of $S$. For connected $S$, this implies all pairs of
vertices in $\Gamma(S)$ are already adjacent to one another.
$\Box$
(ii) ($\Leftarrow$) If $G-S_1 = G\ast S_1$ and $G-S_2 = G\ast S_2$
then $G-S_1-S_2 = G\ast S_1-S_2 = G\ast S_1 \ast S_2$, and $S$ is
the union of its components.
($\Rightarrow$) By definition each $S_i$
is a subset of $S$ and not adjacent to any other connected component of $S$, therefore
$\Gamma(S_i) \subset \Gamma(S)$. If $G-S = G\ast S$ then shorting $S$ in $G$ introduces
no new edges among $\Gamma(S)$, therefore it introduces no new edges among any
subset of $\Gamma(S)$, including $\Gamma(S_i)$. Therefore shorting any part of $S$ does
not introduce new edges among $\Gamma(S_i)$. Therefore $G\ast S_i = G-S_i$.
The argument applies to all the components. $\Box$
(iii) Let $V_1 = V \m S$ and $V_2 = V \cap S$; we will prove
the result for $V_1$ and $V_2$ separately, from which its validity for
$V = V_1 \cup V_2$ follows.
For $V_1$, the result follows immediately by applying the
operation $-V_1$ or $\ast V_1$ to both sides of the equation $G \ast S = G-S$ and using commutation.
For $V_2$, use parts (i) and (ii), as follows. For each connected component $S_i$,
if $\Gamma(S_i)$ is a clique
in a given graph $G$ then $\Gamma(S_i \m V_2)$ is a clique in $G-V_2$
because a subset of a clique is a clique; this is sufficient to
prove the result for $G-V_2$. For $G\ast V_2$ use the fact that the
operation $\ast V_2$ consists in adding edges and then cutting $V_2$.
However if each $\Gamma(S_i)$ is a clique then for $V_2 \subset S$
any edges introduced by the operation $\ast V_2$
must be between members of $S$. It follows that the neighbourhood of each $S_i$ is not
influenced by the addition of those edges and is still a clique. Now re-use the
argument invoked for $G-V_1$. $\Box$

\begin{lemma} \label{lm:connect}
A vertex set $B$ is threatened or supported by $A$ if and only if each
of the connected components of $B$ is respectively threatened or supported
by $A$.
\end{lemma}

Proof: Apply lemma~(\ref{lm:clique})(ii) to vertex set $B$ and the graph $G-A$ or $G\ast A$,
for threat or support respectively. $\Box$

{\bf Definition} If for a vertex $w$ and vertex set $U$, $\Gamma(U) \subset \Gamma[w]$
(i.e. the closed neighbourhood of $w$ contains the 
neighbourhood of $U$), then $w$ is said to {\em surround} $U$. For
example, $w$ surrounds a single vertex $v$ if and only if
$\Gamma(v) \subset \Gamma[w]$. Also, $w$ surrounds the 
set $\Gamma[v]$ if and only 
if $\Gamma( \Gamma(v) \setminus \Gamma(w) ) \subset \Gamma[v] \cup \Gamma(w)$.

\begin{theorem} \label{th:threat_support}
(i) Set $A$ threatens a connected set $B$ if and only if $\Gamma(B) \m A$ is a
clique.\\
(ii) Set $A$ supports a connected set $B$ if and only if
$A$ is adjacent to all neighbours of $B$ that do not surround $B$.
(This property can be expressed $\Gamma(B) \m g(B) \subset
\Gamma(A)$, where $g(B)$ is the set of all neighbours of $B$ that
are adjacent to all other neighbours of $B$).
\end{theorem}

Proof: (i) ($\Leftarrow$) After deleting $A$, $B$ satisfies the
clique condition of lemma \ref{lm:clique}(i), hence $G-A\ast B=G-A-B$.
($\Rightarrow$) If $G-A\ast B=G-A-B$ then all edges introduced
by shorting $B$ must be incident on $A$.
Therefore $\Gamma(B_i) \m A$ is a clique for each connected component
$B_i$ of $B$. Since $B$ is connected there is only one such component
$B_1 = B$. $\Box$ \\
(ii) ($\Leftarrow$) $A$ is adjacent to all those neighbours of $B$
that are not adjacent to all other neighbours of $B$. It follows
that shorting $A$ will convert the neighbours of $B$ into a
clique, and the result follows from lemma \ref{lm:clique}.
($\Rightarrow$) The property is $(G \ast A)-B =(G \ast A)\ast B$
with connected $B$, therefore in graph $G \ast A$, the neighbours of $B$ form a
clique (lemma \ref{lm:clique}). It follows that, in $G$, shorting $A$ must cause
all non-adjacent neighbours of $B$ to become adjacent, from which it follows that
$A$ must be adjacent to all neighbours of $B$ that are not
adjacent to all other neighbours of $B$ (and $A$ is allowed to
have further neighbours). $\Box$.

\begin{figure}
\centerline{\resizebox{!}{0.1\textwidth}{\includegraphics{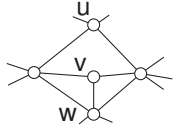}}}
\caption{Example of `surround' and `support'.}
\label{f.surround}
\end{figure}

The notions of `surrounding' and `supporting' are related but different.
If $w$ surrounds $v$ then $w$ supports $v$, but a vertex can support another
without necessarily surrounding it. In figure \ref{f.surround},
$u$ supports $v$ but does not surround it; $w$ surrounds (and therefore also
supports) $v$.

{\bf Examples}
Any vertex set is threatened by its neighbourhood, and by its
neighbourhood minus a clique.
A vertex of degree 2 is threatened by each of its
neighbours, and supported by any vertex to which it has more than one
2-walk.

Theorem \ref{th:threat_support} enables supports and threats
to be identified easily, especially single-vertex
supports or threats. For a single vertex $b$,
one way to identify the set $g(b)$ is to use the fact that it is
the intersection of $\Gamma(b)$ with the set of vertices having
$(d-1)$ 2-walks to $b$, where $d$ is the degree of vertex $b$.
In fact one can avoid the need to calculate $g(b)$, because when
the theorem applies, the edges from $b$ to $g(b)$ will all be dead
by theorem \ref{th:vwdead}, so $g(b)$ is an empty set after
deletion of dead edges (defined below) and then the condition of
theorem \ref{th:threat_support}(ii) simplifies to $\Gamma(b) \subset
\Gamma(a)$. This can easily be checked using an adjacency matrix
or adjacency lists.

To discover threats to a single vertex, one can use 2-walks and
the triangle number:

\begin{theorem} \label{th:triangle}
For a given vertex $v$, let $x_{n}(v) = T(v) - (d-n)(d-n-1)/2$, where $T$ is
the triangle number and $d$ the degree of $v$. $v$ is threatened by
any set of $n$ of its neighbours having a total of $x_n(v)+m$ 2-walks
to $v$, where $m$ is the number of edges between members of the set.
\end{theorem}

{\bf Corollary}: A single neighbour having $(T(v) - (d-1)(d-2)/2)$ 2-walks to $v$
threatens $v$.

Proof: deletion of such a set $S$ would lower the
degree of $v$ by $n$ and the triangle number by $x$. The new triangle number
shows that $v$ is now only adjacent to a clique $K_{d-n}$\footnote{In any
graph, if a vertex of degree $d$ has triangle
number $d(d-1)/2$ then that vertex is only adjacent to a
clique $K_{d}$, because to get this many
triangles all neighbours of the vertex must be adjacent to each other.}.
Therefore, using lemma \ref{lm:clique}(i),
$(G-S) \ast v = (G-S)-v$ which is the condition defining a threat.  $\Box$

Having gone to the trouble of identifying a threat or a support, the
issue arises, will the property be preserved during further play, or must
it be re-calculated for every new position?

\begin{lemma}  \label{lm:robust}
(i) If $A$ threatens $B$ in $G$ then
$A$ threatens $(B\m \{b\})$ in $G-b$ and in $G \ast b$, $\forall
b \notin A$.\\
(ii) If $A$ threatens $B$ in $G$ then $(A-a)$ threatens $B$ in $G-a$, $\forall a \in A$.\\
(iii) If $A$ supports $B$ in $G$ then
$A$ supports $(B\m \{b\})$ in $G-b$ and in $G \ast b$, $\forall b
\notin A$.\\
(iv) If $A$ supports $B$ in $G$ then
$(A-a)$ supports $B$ in $G\ast a$, $\forall a \in A$.
\end{lemma}

{\bf Corollary}
For distinct sets $A,B_1,B_2$,
if $A$ threatens (supports) $B_1 \cup B_2$ in $G$ then $A \cup B_1$ threatens
(supports) $B_2$ in $G$.

Proof: (i),(iii) apply lemma \ref{lm:clique}(iii) to the set
$B$ in the graphs $G-A$ and $G\ast A$. $\Box$ (ii) Note that $G-A \equiv
(G-a)-(A-a)$ and apply it to both sides of $G-A\ast B = G-A-B$;
the result follows immediately using commutation. $\Box$ (iv) use
$G\ast A \equiv (G\ast a)\ast (A-a)$ and proceed as in (ii).
$\Box$ The corollary can be obtained from (i) and it also follows
immediately from lemma \ref{lm:clique}(iii) applied to
the graph $G-A$ or $G\ast A$. $\Box$

Parts (i) and (iii) of the lemma show that a threat or a
support can only be removed by playing directly on the threatening
or supporting set. This is in contrast to domination,
since if a vertex $v$ $P$-dominates a set $S$, a play by $\neg P$ on
some other vertex in $S$ may remove the domination.
In the context of computational analysis of a game the robustness
is useful, because the computational effort invested in
identifying a threat or support does not become valueless
as soon as a new move is made. Parts (ii) and (iv)
show that threat and support also survive play by one player in
the threatening or supporting set. This again helps the computational
problem, and is used in the proof of theorem \ref{th:threaten3}. The corollary
is used in theorem \ref{th:pairset}.

% that if $A$ threatens $B$ and $B$ threatens $C$, then $A$ threatens
%$B \cup C$:
%\begin{eqnarray}
%& G-A \ast B &=\, G-A-B, \;\;\ G-B\ast C = G-B-C, \\
%\Rightarrow & G-A \ast B \ast C = G-A-B-C
%\end{eqnarray}

Next, observe that the action of `removing' a threat by shorting the
threatening set $A$ really transfers the threat to the
neighbourhood of $A$:
\be
\lefteqn{G - A \ast B = G - A - B}  \nonumber \\
\lefteqn{\Rightarrow \;\;\; G - (\Gamma^G(A) \m B) \ast A \ast B =}  \nonumber \\
&& G - (\Gamma^G(A) \m B) - A - B
\label{nthreat}
\ee
where the superscript $G$ is to
indicate that the neighbourhood in question is that in the
original graph $G$, not the one obtained after shorting or
deleting $A$ and $B$. In words the result is:
if $A$ threatens $B$ then $\Gamma(A) \m B$ threatens
$A \cup B$. To obtain the result, use the property of threat, that
new edges introduced by shorting $B$ are all incident on $A$.
Therefore shorting $B$ followed by cutting $(\Gamma^G(A) \m B)$
is equivalent to (results in the same graph as) cutting $B$ and
$(\Gamma^G(A) \m B)$, i.e. cutting $\Gamma^G(A)$, and of course
$\Gamma^G(A)$ threatens $A$.

This permits a straightforward proof of the following.

\begin{lemma}
A vertex of degree 3 that threatens each of two of its neighbours is in a set
cut-dominated by the other neighbour.
\end{lemma}

\begin{theorem}  \label{th:threaten3}
A vertex of degree 3 that threatens all its neighbours
is in a lost set (and therefore so are its neighbours).
\end{theorem}

Proof. First we consider the lemma. Let $v$ be the vertex of degree 3, let $a,b$
be the threatened neighbours and let $w$ be the non-threatened neighbour.
After shorting $v$, the threat that it offered to $a$ and $b$
is transferred as in eq. (\ref{nthreat}), such that $a$ is now threatened
by the set $\{b,w\}$ and $b$ is threatened by the set $\{a,w\}$. 
After cutting $w$,
the resulting situation is one in which each of the remaining vertices
threatens the other (lemma \ref{lm:robust}(ii)), so both are lost
(corollary to theorem \ref{th:domts}).
It follows that $w$ cut-dominates $\{v,a,b,w\}$.

To prove the theorem, we show that \cut~has a second player winning
strategy for ${\cal M}_{G,S}$ where $S=\{v,a,b,w\}$, as
follows. If \short~does not short $v$, then \cut~deletes it and
wins (since it threatens all the others, the graph is now
$G-S$). Otherwise, after shorting $v$, \cut~may delete any vertex
and thus leave a mutually threatening pair, hence a lost pair. $\Box$

Note that theorem \ref{th:threaten3} allows the reduction of a 4-vertex
subgraph that is not reducible by theorem \ref{th:doubledom} alone.

There is no partner result to eq. (\ref{nthreat}) following the
removal of a support. In this case there does not necessarily
remain any supporting set of vertices (one may easily show that a
vertex has no supporting set if and only if it is an articulation
vertex (cut-vertex) of the graph). When there is no single supporting
vertex, but a larger supporting set exists,
there is no simple rule for finding a small or minimal such set.

Some further captured or lost sub-graphs not reducible by
theorem \ref{th:doubledom} are covered by the following theorem.

\begin{theorem} \label{th:pairset}
For any set of $P$-dominated vertex pairs $\{a_i,b_i\}$, let $A$
be the set of dominating vertices and $B$ the set of dominated
vertices (i.e. $A = \Cup_i \{a_i\} $, $B  = \Cup_i \{b_i\}$ where $a_i$
$P$-dominates $\{a_i,b_i\}$). If for such a set of pairs, the
pairs are cut-dominated and $B$ threatens $A$, then $A \cup B$ is
lost; if the pairs are short-dominated and $B$ supports $A$, then
$A \cup B$ is captured.
\end{theorem}

\begin{figure}[ht!]
\centerline{\resizebox{!}{0.1\textwidth}{\includegraphics{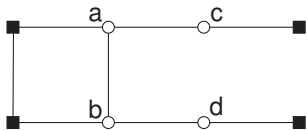}}}
\caption{Example multi-Shannon game that is not
reducible by mutually dominating pairs, but can be solved
by pair-sets (theorem \protect\ref{th:pairset}). In the game shown,
$a$ threatens $c$ and $b$ threatens $d$ and $\{c,d\}$ threatens $\{a,b\}$
hence there is second-player winning strategy for \cut, by
theorem \ref{th:pairset}, so the set $\{a,b,c,d\}$ is lost.}
\label{f.eg_pairset}
\end{figure}

Proof: For brevity we only
present the proof for $P$=\cut; the proof for $P$=\short\ is
similar. The 2nd-player winning strategy for \cut~is as follows.
Each move by \short~must be on one of the pairs $a_i, b_i$;
\cut~replies by cutting the other vertex of the pair. Let $G_n$
be the reduced graph at the $n$'th pair of moves.
On each occasion that a vertex in $A$ is shorted, the graph develops to
$G_{n+1} = G_{n} \ast a_i - b_i$. On each occasion that a vertex in
$B$ is shorted, the graph develops to $G_{n+1} = G_n - a_i \ast b_i
= G_n - a_i - b_i$ by the definition of {\em threat}
(see theorem \ref{th:domts}). Hence after all
vertices are coloured, the reduced graph is the same as one in which
all of $B$ is cut. Let $A_1 \subset A$ be the part (possibly empty)
of $A$ that was shorted, then obviously $A_0 \equiv A \m A_1$ was cut,
so the final reduced graph is $G \ast A_1 - A_0 - B$.
By the corollary to lemma \ref{lm:robust}, $B \cup A_0$ threatens $A_1$, therefore
this reduced graph is equal to $G-A-B$, therefore the multi-Shannon
game is won by \cut. $\Box$ An example is shown in figure \ref{f.eg_pairset}.

%Proof (by induction): For $m=1$ the statement is true by theorems
%\ref{th:doubledom}, \ref{th:domts}.
%It remains to prove the theorem for sets of size $m$
%assuming it is valid for sets of size $m-1$. For clarity we
%present the proof for $P$=\cut; the proof for $P$=\short\ is
%similar. The 2nd-player winning strategy for \cut~is as follows.
%An opening move by \short~must be on one of the pairs $a_i, b_i$;
%\cut~replies by cutting the other vertex of the pair. If $a_i$ was
%cut, then the graph is now $G-a_i\ast b_i = G-a_i-b_i$ by
%the definition of {\em threat}.
%If $b_i$ was cut, then the graph is now $G \ast a_i -
%b_i$. In either case, the graph is equal to one obtained from
%$(G-b_i)$ by shorting or cutting a vertex in $A$. Using lemma
%\ref{lm:robust}(ii), $(B-b_i)$ threatens $A$ in $G-b_i$, and
%therefore $(B-b_i)$ threatens $(A-a_i)$ in $G-b_i$ ?????. It
%follows that both the cases $(G-b_i) - a_i$ and $(G-b_i) \ast a_i$
%realise the condition of the theorem for sets $(A-a_i), (B-b_i)$
%of size $m-1$. Assuming the theorem for $m-1$, we can now cut all
%of $A-a_i$ and $B-b_i$ without influencing the outcome of the
%multi-Shannon game on $A \cup B$. What remains is then either $G -
%A - B$ or $G - (A - a_i) \ast a_i - B$. Using the fact that $B$
%threatens $A$ in $G$, the latter is the same as $G-A-B$. It
%follows that the game outcome is the same as if all of $A$ and $B$
%were deleted, therefore it is won by \cut, for any opening move by
%\short. $\Box$

Finally, we note a \cut-winning condition for the Shannon game
that is easy to check if threats are in any case being examined.
(Simple 2nd-player \short-winning conditions are well known, such as
the presence of more than one 2-walk between terminals.
By considering the block graph one can derive further conditions for
a \cut\ win, but to check this is relatively costly in computational terms.)

\begin{theorem}
In the Shannon game, if the neighbours of a terminal $t$ are all threatened by distinct vertices,
then \cut~has a 2nd-player winning strategy. The strategy is: if \short~shorts
a vertex in $\Gamma(t)$, then cut the corresponding threat vertex; if \short~shorts
a vertex threatening $v \in \Gamma(t)$, then cut $v$, otherwise cut any
remaining vertex in $\Gamma(t)$.
\end{theorem}

{\bf Corollary} If all the neighbours of a terminal except one
are threatened by distinct vertices, then \cut\ has a 1st-player winning strategy, namely
to cut the non-threatened neighbour.

{\bf Example} If a terminal $t$ has neighbours all of degree 2, and no
vertex has more than one 2-walk to $t$, then \cut~is the winner.

\subsection{Discussion}

The theorems concerning threat and support  simplify the analysis of 
the Shannon game, and the computational task of identifying
lost or captured sets.

We described the multi-Shannon game in section \ref{s:multi}
because this is the right way to understand capture and loss in general. However,
we don't need that concept to understand support and threat. If we ignore
the multi-Shannon game we can still prove the essential property
of threat pairs and supporting pairs, namely that $W(G,x) =
W(G\ast a \ast b,x)$ if $a$ and $b$ support each another, and
$W(G,x) = W(G - a - b,x)$ if $a$ and $b$ threaten each other (for single
vertices $a$, $b$). For,
if $a$ and $b$ support each other then $G - a \ast b = G \ast a - b
= G \ast a \ast b$. Therefore if $W(G,x)=$\cut~(for some value,
not necessarily all values, of $x$) then there is a
winning strategy whose cutset does not require $a$ or $b$, hence they are
spectators. If the winning cutset did require one of $a$ or $b$, then at some stage
\cut~must cut it, but then \short~can reply such that the position
is equivalent, in the colouring model, to one in which both are coloured black
(\short's colour) and we have a contradiction. If $W(G,x)=$\short~then
obviously $W(G\ast a\ast b,x)=$\short~since shorting vertices as
a ``free move'' can never be disadvantageous to \short~(this is obvious
in the graph colouring model). The same
argument with appropriate adjustments proves the corresponding result
for a threat pair. $\Box$

\begin{figure}[ht!]
\centerline{\resizebox{!}{0.2\textwidth}{\includegraphics{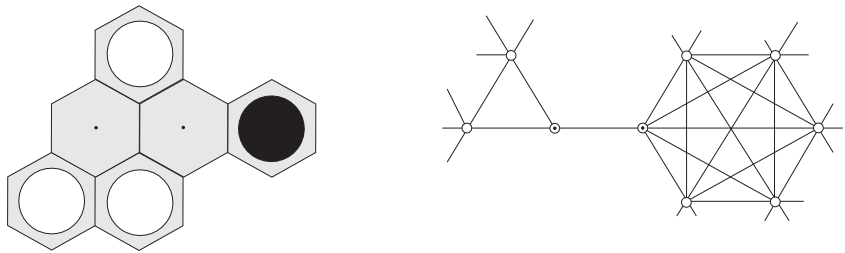}}}
\caption{Example of a threat-pair arising in Hex. The left diagram shows part of
a Hex board, with \short~playing black.
The right diagram shows the equivalent graph, supposing the cells not
shown on the left are empty. Theorem \protect\ref{th:triangle}
can be used to show that the dotted cells (vertices) threaten each other.
If the surrounding cells are not empty, the effect on the reduced graph
is merely to grow or shrink the cliques, see lemma \ref{lm:robust}.}
\label{fig:eghexloss}
\end{figure}

The graph properties identified in
theorems \ref{th:threat_support} and \ref{th:triangle} suffice to
discover (through the above argument, or by using
theorems \ref{th:domts}, \ref{th:doubledom}) all the examples of
lost pairs presented in \cite{06Bj} (they are called `cut-captured' there),
and all the examples of captured pairs
presented in \cite{06Bj}, without the need for any further pattern search.
An illustrative example is shown in figure \ref{fig:eghexloss}.

\begin{figure}[ht!]
\centerline{\resizebox{!}{0.15\textwidth}{\includegraphics{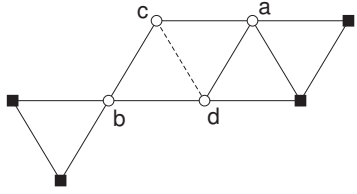}}}
\caption{Example multi-Shannon game that can arise in Hex, and that
can be solved by pair-sets (theorem \protect\ref{th:pairset}).
First note that $a$ threatens $c$, and observe from theorem
\protect\ref{th:vwdead} that the edge $cd$ is dead. After deleting it,
$b$ threatens $d$.
Also, $\{c,d\}$ threatens $\{a,b\}$ (since $\Gamma(a) \m \{c,d\}$
is a clique and $\Gamma(b) \m \{c,d\}$ is a clique, lemma \protect\ref{lm:connect}
and theorem \protect\ref{th:threat_support}). Therefore the set of
pairs $\{ \{a,c\}, \{b,d\} \}$ realise the condition of
theorem \ref{th:pairset} and therefore $\{a,b,c,d\}$ is lost.}
\label{fig:lmgame_abcd}
\end{figure}

\begin{figure}[ht!]
\centerline{\resizebox{!}{0.15\textwidth}{\includegraphics{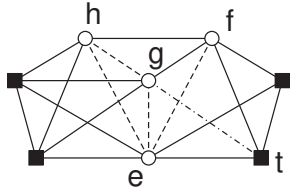}}}
\caption{A multi-Shannon game that can arise in Hex and that is not
reducible by mutually dominating pairs, but can be solved by using
{\sc and}-, {\sc or}-rules for multiple connections, see appendix,
or by using supporting sets. There is a second-player winning strategy
for \short, hence vertices $\{e,f,g,h\}$ are `captured'. Proof:
$e$ supports each of $f$, $g$, $h$ (theorem \protect\ref{th:threat_support})
so each of them is locally losing for \cut.
$f$ short-dominates $\{f,g,h\}$ since after shorting $f$, $g$ and $h$ support
each other and are therefore captured (theorems \protect\ref{th:domts},
\protect\ref{th:doubledom}).
Also, $\{f,g,h\}$ supports $e$. Therefore we have that \short\ dominates
a set that supports $e$, hence $e$ is locally losing for \cut. Thus all moves are locally
losing for \cut, therefore \short\ has a second player winning strategy.
(The edges shown dashed are dead, and this
can be discovered using theorem \protect\ref{th:vwdead}, but we didn't need to invoke
that here. The edge $gt$ (shown dash-dotted) is not dead but it is a spectator. If its
spectator status could be discovered conveniently (without first solving the game) then
it could be deleted, and then the game could be reduced using support pairs.)}
\label{fig:mgame_efgh}
\end{figure}

Bj\"ornsson {\em et al.} also presented some examples of won multi-Shannon
games that arise in Hex but are not reducible by mutually dominating pairs.
All of these are covered by theorems \ref{th:threat_support},
\ref{th:threaten3} and \ref{th:pairset}.
Figures \ref{fig:lmgame_abcd} and \ref{fig:mgame_efgh} present two examples.
These arise from two
dual graphs representing the
same hex position; figure \ref{fig:lmgame_abcd}
is a \cut-win; figure \ref{fig:mgame_efgh} is a \short-win. The captions give
the proofs.

In a computational setting, one can track both
threat and support by maintaining a (representation
of a) directed graph for each. This allows threats
and supports to be detected incrementally and locally as the game proceeds (at any
moment during the course of a game, one only
has to examine the neighbourhood of the part of the graph which just changed).
Theorems \ref{th:threat_support} and \ref{th:triangle} enable
threats and supports to be detected rapidly.
Once detected they are stable (lemma \ref{lm:robust}). Suppose for example
that at some stage it is found that $a$ threatens $b$.
Such a threat remains until such time as $a$ is shorted, so it does not need
to be reconsidered as the game proceeds. Thereafter neither player should
play at $b$, until such time as $a$ is shorted. This is because if \cut~wants
to delete $b$, he may as well delete $a$ and then $b$ dies, and \short~can be
sure that a play at $b$ is useless to him, since \cut~could, if he chose, reply
at $a$ which would return a position equivalent to the one \short~was considering,
but with two vertices deleted.

If subsequently a position develops in which 
$b$ also threatens $a$, then both can be deleted
immediately. Similar considerations apply to a supported vertex,
with the roles reversed.

\section{Dead edges and vertices} \label{s.dead}

%We will say that two
%graphs $G,G'$ (with their terminal designations) are `equivalent' if
%$W(G,x) = W(G',x) \forall x$, i.e. their win-values are the same when the same
%player is to move. Such equivalence is indicated $G \cong G'$.

If $W(G-v,x) = W(G \ast v,x) \forall x$ then we say vertex $v$ is
a {\em spectator}: no matter who is to play, the game outcome is
not influenced by deleting or shorting $v$.

A stronger property is discussed in \cite{06Bj}. In the graph
colouring model, suppose at a given game position we colour all
the remaining uncoloured vertices in some arbitrary black and
white pattern, without respect to taking turns. Such a complete
colouring $c$ produces some final position ${\cal G}_c$ with a
definite winner.

{\bf Definition} In a given graph $(G,\{t_1,t_2\})$ a non-terminal
vertex $v$ is {\em dead} if, for every complete colouring ${\cal
G}_c$ obtainable from $G$, the winner of ${\cal G}_c$ is
independent of the colour of $v$.

{\bf Definition} In a given graph $(G=(V,E),\{t_1,t_2\})$ an edge
$e \in E$ is {\em dead} if, for every complete colouring ${\cal
G}_c$ obtainable from $G$, the winner of ${\cal G}_c$ does not
change if $e$ is deleted. An edge $e \notin E$ is dead if it is
dead in $(G+e,\{t_1,t_2\})$.

The dead vertex idea was discussed in \cite{06Bj}; here we extend
the concept to a dead edge. This will allow us to introduce some
new useful results.

%We will say that two graphs $G,H$ are `equivalent' if
%they are the same {\em after deleting dead edges}, or isomorphic
%after deleting dead edges by an isomorphism which maps
%the terminal set of $G$ to the terminal set of the other $H$.
%Such equivalence is indicated $G \cong H$.

\begin{theorem} \label{th:dead_ev}
A vertex is dead if and only if all edges incident on it (and
present in the graph) are dead.
\end{theorem}
Proof: Consider the final position after some complete colouring
obtainable from the given initial position. ($\Leftarrow$) If all
the edges incident on $v$ are dead, then we can delete them
without changing the winner, and then $v$ is isolated so it is
obvious that its colour does not influence any outcome.
($\Rightarrow$) If the vertex is dead, then if the winner is
\short~the vertex is not needed so its edges are irrelevant, and
if the winner is \cut~it is still \cut~if any edge is deleted. In
either case the outcome is unchanged if any of the incident edges
are deleted.  $\Box$

{\bf Corollary} Adding or deleting dead edges in a game $G$ does
not influence whether or not any vertex is dead in $G$.

Proof: If a vertex is dead, then after adding or deleting dead
edges it still satisfies the the sufficient condition of the
theorem. If a vertex is not dead, then adding or deleting dead
edges does change the fact that it has a non-dead edge so does not
satisfy the necessary condition of the theorem.

In a given game position, the set of dead vertices can include
some that are already coloured by one or other player (but
changing their colour will not affect the outcome). If a vertex is
dead then it is a spectator.

Dead vertices are in general hard to identify: recognising dead
vertices in the Shannon game is NP-complete \cite{06Bj}. In this
section we present some theorems that allow some dead vertices and
dead edges to be identified easily, using readily calculated
properties of the graph. The first case, a subgraph with a clique
cutset (theorem \ref{th:vdead}), is already known, see for example
\cite{06Hayward}, it is included here for completeness and because
the proof is instructive.

\begin{theorem}  \label{th:vdead}
A terminal-free subgraph with a clique cutset is dead (i.e. all its edges and
vertices are dead).
\end{theorem}

Proof: The vertices of the subgraph are a subset of some vertex set $V$ whose
neighbourhood is a clique. For such a set
$G \ast V = G - V$. We will prove that all the members of $V$ are dead.
Pick a vertex $v \in V$. First consider some complete colouring $\cal C$
in which $v$ is black (\short's colour). Then, if the winner is \cut\ it is obviously
still \cut\ if $v$ is changed to white. If the winner is \short\ then it is
still \short\ if all the rest of $V$ is coloured black (since adding edges cannot
be bad for \short). Since the reduced graph $G\ast V$  is the same as $G-V$, the winner
must still be \short\ if all of $V$ is white. Adding edges to
the reduced graph $G-V$ cannot be bad for \short, so the winner remains \short\ if
we return to the colouring $\cal C$ with
just the colour of $v$ changed from black to white. This completes the proof
for complete colourings
in which $v$ is black. A similar argument, adjusted appropriately,
applies for complete colourings in which $v$ is white. $\Box$

{\bf Example} 
Non-terminal simplicial\footnote{A vertex is called {\em simplicial} if its
neighbourhood is a clique; or to be precise, if the subgraph induced by its
neighbourhood is a complete graph.} vertices are dead.
(For example, pendant vertices are simplicial.)

\begin{figure}[ht!]
\centerline{\resizebox{!}{0.125\textwidth}{\includegraphics{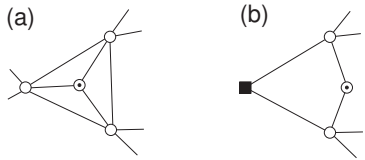}}}
\caption{Examples of dead vertices. A vertex shown with a dot is
dead. A filled square vertex is a terminal in a Shannon game. (a)
illustrates theorem \protect\ref{th:vdead}; (b) illustrates
theorem \protect\ref{th:dead_tneighbour} and is an example where
$G \ast v \ne G-v$. } \label{fig:dead_vertex}
\end{figure}

Clique cutsets can be discovered by an $O(n^3)$
algorithm due to Whitesides \cite{81Whitesides,77Gavril}. There can be
dead vertices which do not satisfy the condition of theorem \ref{th:vdead};
figure \ref{fig:dead_vertex}(b) shows an example of a dead vertex
where $G \ast v \ne G-v$ (see theorem \ref{th:dead_tneighbour}).

{\bf Definition} An edge between a vertex $v$ and a surrounding 
vertex (i.e. one surrounding a set containing $v$) is called {\em transverse}.

\begin{theorem}  \label{th:vwdead}
If $w$ surrounds a terminal-free subgraph containing
$v$, then $vw$ is dead. More succinctly, transverse edges to
a terminal-free subgraph are dead. 
\end{theorem}

{\bf Corollary I} If $w$ surrounds non-terminal $v$ then $vw$ is dead. 
(This condition is equivalent to
\[
A^2(v,w) = d(v) - A(v,w)
\]
where $A$ is the adjacency matrix and $d$ is degree.)

{\bf Corollary II} 
If $\Gamma[v]$ is terminal-free and
\be
A^3(v,w) = \sum_{i \ne w} A^2(v,i).     \label{second_surround}
\ee
then all edges between $w$ and $\Gamma[v]$ are dead. 

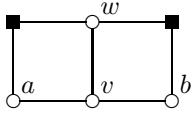
\begin{figure}
\begin{picture}(100,40)(-20,0)
\multiput(-2.5,27.5)(60,0){2}{\rule{5pt}{5pt}}
\put(30,30){\circle{5}}
\multiput(0,0)(30,0){3}{\circle{5}}
\multiput(2.5,0)(30,0){2}{\line(1,0){25}}
\multiput(2.5,30)(30,0){2}{\line(1,0){25}}
\multiput(0,2.5)(30,0){3}{\line(0,1){25}}
\put(33,33){$w$}  \put(3,3){$a$}  \put(33,3){$v$}   \put(63,3){$b$}
\end{picture}
\caption{Edge $vw$ is dead (by corollary II to theorem \ref{th:vwdead}); after
deleting it, one may discover by theorem \ref{th:doubledom} (threat pairs)
that the pairs $\{a,v\}$ and $\{v,b\}$ are lost.}
\label{f.dead_edge_vw}
\end{figure}

The theorem is based on the idea that if a winning path passing
from $w$ to $v$ must continue to a neighbour of $w$ before it can
reach a terminal, then the part through $v$ must be irrelevant.
(This is loosely reminiscent of, but clearly different from, lemma 1 of
Hendersen {\em et al.} \cite{Henderson:2009:SH:1661445.1661526}).
The corollaries give two examples of the condition of the
theorem that are easy to identify in practice. An example
of corollary I is shown in figure
\ref{fig:dominate}, where vertex $c$ surrounds vertex $a$. In corollary II, all walks
from $v$ either remain in $\Gamma[v]$ or pass through $\Gamma[w]$; figure
\ref{f.dead_edge_vw} shows an example.

Proof: invoke the graph-colouring model, and let $e=vw$ be the
edge. For each complete colouring invoked in the definition of
dead edges, the final position is either won or lost for Short. If
it is lost then obviously the deletion of $e$ will not change the
win-value. If it is won then $e$ is only alive if the existence of
a winning path requires it. However, for every short-winning path
$\left< \ldots, w,v, P, a, \ldots \right>$ containing $e$, where
$P$ is a sub-path (possibly empty), there is a winning sub-path
not containing $e$: $\left< \ldots, w, a, \ldots \right>$ (as long
as neither $v$ nor the vertices of $P$ are terminal). $P$ refers
to that part of the path after $v$ which remains on the subgraph
surrounded by $w$. Eventually the winning
path must emerge from this vertex set, and when it does the next
vertex $a$ is a neighbour of $w$, by the condition of the theorem.
$\Box$

Proof of eqn (\ref{second_surround}). Corollary II presents a
sufficient but not necessary condition that $w$ surrounds $\Gamma[v]$.
Consider all $2-$walks from $v$. 
Let $N_1$, $N_2$, $N_3$ be the number that end on a neighbour of $w$,
on $w$ itself, and on neither, respectively. Then (since the $i$'th power of
the adjacency matrix gives the number of $i$-walks),
\[
\sum_i A^2(v,i) = N_1 + N_2 + N_3
\]
and
\[
A^3(v,w) = N_1, \;\;\;\;\; A^2(v,w) = N_2.
\]
If $N_3=0$ then $w$ surrounds $\Gamma[v]$, since then all walks from
$v$ of length 2 or more start with a part that passes through $\Gamma[w]$.
(This is a sufficient but not necessary condition; it ignores the possibility
that there may be a triangle including $v$ and a non-neighbour of $w$, and also
the possibility of two-walks from $v$ that pass through $\Gamma(w)$ 
to a vertex other than $w$).
Substituting $N_3=0$ into the above set of equations gives
(\ref{second_surround}).

The discovery of even a single dead edge can considerably simplify
analysis of a given Shannon game. For example, figure \ref{f.dead_edge_vw}
shows a case
where after deleting one dead edge, three vertices can also be
eliminated because it becomes obvious that they form a lost set. In the context
of a larger game, this in turn reduces the amount of further information
needed to be handled by a solution algorithm, such as the set of
vertex sets forming weak and strong `links' or `virtual connections' \cite{06Hayward}.

We conclude this section with two theorems which stem from the special role 
played by terminals in the Shannon game.

\begin{theorem} \label{th:dead_tedge}
An edge between neighbours of a given terminal is dead.
\end{theorem}

Proof: An edge $uv$ is only alive if the existence of a
\short-coloured path between the terminals requires it in
some complete colouring. Let $t$ be the terminal mentioned in the
theorem. All paths sought by \short\ start or end on $t$. 
In the condition of the theorem, both $u$ and $v$ are in
$\Gamma(t)$, and therefore any path $\left<t,\ldots u,v,
\ldots\right>$ has a subpath $\left<t,v, \ldots\right>$, and any
path $\left<t,\ldots v,u, \ldots\right>$ has a subpath $\left<t,u,
\ldots\right>$, therefore edge $uv$ is not required by \short~in
order to form any path to or from $t$. Therefore
$uv$ is not required by \short\ in any complete colouring, so is
dead. $\Box$

\begin{theorem} \label{th:dead_tneighbour}
A terminal-free vertex set surrounded by a terminal is dead.
\end{theorem}

(Figure \ref{fig:dead_vertex}b shows an example.)
We offer two proofs. Proof 1: The closed or open
neighbourhood of the terminal can be converted into a
clique by adding edges that are dead by theorem
\ref{th:dead_tedge}. After this is done the vertex set in question
has a clique
cutset so is dead by theorem \ref{th:vdead}, and therefore it is
dead in the original graph by the corollary to theorem
\ref{th:dead_ev}. $\Box$

Proof 2: Let $U$ be the vertex set in question.
For a vertex to be alive, it is necessary that a path
sought by \short~requires it in some complete colouring. However,
all paths sought by \short~start
or finish on $t$, and 
in the conditions of the theorem,
any path $\left<t,\ldots P,a,
\ldots\right>$, where $P$ is on $U$, has a subpath $\left<t,a,
\ldots\right>$, since any non-$t$ neighbour of $U$ is a neighbour
of $t$. Therefore the vertices of $U$ are dead. $\Box$

\section{Statistics on small graphs}

The usefulness or otherwise of the theorems presented in this paper, for the purpose
of analyzing the Shannon game, can be roughly assessed by exploring how
often the conditions of the various theorems arise in Shannon games in
general. To this end, we examined all graphs on up to 10 vertices, 
and counted those in which the conditions of various of the theorems arise.
Table \ref{t.count} shows the results. The main thing to notice is that 
the numbers in columns 3--6 of the table are small compared to the number
of connected graphs (column 2). This shows that the great
majority of graphs satisfy one or more of the conditions, so the theorems
will be useful in practice to simplify most games. 

%The single most useful theorem is arguably theorem \ref{th:vwdead}.

\begin{table}
\begin{tabular}{r|rr|rrrr}
size & graphs & S-free & $\perp$-free & $\ge 2 \,\triangle$-free & 4 \& 5\\
\hline
 1  &        1  &        0  &        1  &        0  &        0  \\
 2  &        1  &        0  &        0  &        1  &        0  \\
 3  &        2  &        0  &        0  &        1  &        0  \\
 4  &        6  &        1  &        1  &        3  &        1  \\
 5  &       21  &        4  &        2  &       10  &        2  \\
 6  &      112  &       24  &        9  &       52  &        7  \\
 7  &      853  &      191  &       46  &      363  &       34  \\
 8  &    11117  &     3094  &      507  &     4022  &      327  \\
 9  &   261080  &    95204  &    11800  &    72594  &     5897  \\
10  & 11716571  &  5561965  &   626586  &  2276219  &   213064  \\
\hline		  
\end{tabular}
\caption{Numbers of graphs with various properties. In all cases only connected graphs
are counted. The second column shows the total number of connected graphs on the given
number of vertices. `S-free' means simplicial-free (i.e. no simplicial vertices). 
`$\perp$-free' means there are no transverse edges (i.e. those which are
dead by corollary I to theorem \protect\ref{th:vwdead} if the surrounded vertex
is non-terminal). The 5th column
counts graphs containing at least 2 triangle-free vertices---only Shannon games on
such graphs are free of edges that are dead by 
theorem \protect\ref{th:dead_tedge}. The
last column counts graphs satisfying conditions of both columns 4 and 5.}
\label{t.count}
\end{table}

\section{Conclusion}

The new ideas presented in this paper are of three kinds.
We presented new methods to identify dead edges and dead vertices:
theorems \ref{th:vwdead},
\ref{th:dead_tedge}, \ref{th:dead_tneighbour}.
Secondly, we presented methods to solve some
multi-Shannon games (theorems \ref{th:threaten3}, \ref{th:pairset}) that are
not easily reducible. 
Thirdly, we introduced the idea of `threat' and `support' (in our terminology).
These are in general stronger conditions than `domination' in the terminology of
\cite{06Bj}. This makes them rarer but easier to identify. In the case of 
a two-vertex set, threat or support is equivalent to domination
(theorem \ref{th:domts}). Therefore to identify mutually dominating pairs
(which are consequently either lost or captured) it is sufficient to
identify mutually threatening or supporting pairs, and we have presented, in
theorem \ref{th:triangle} and corollary I to theorem 
\ref{th:vwdead}, easily-computed properties that can be used to identify
threats and supports without the need for pattern-searching algorithms.

This work was supported by Oxford University.

%\begin{theorem}  \label{support3}
%A vertex of degree 3 which supports three other vertices is captured
%(and therefore so are the supported vertices) if the supported vertices
%form a triangle.
%\end{theorem}

%Proof. Let $v$ be the vertex.
%\short~has a second player winning
%strategy for ${\cal M}_{G,S}$ where $S=\{v,\Gamma(v)\}$, as follows.
%If \cut~does not delete $v$, then \short~shorts it and wins. If \cut~deletes
%$v$, then argue as follows. Before deletion, $v$ has degree 3 and supports
%the whole set. It follows that the multi-Shannon game has at most 3 terminals
%that are non adjacent

\section*{Appendix: {\sc or} and {\sc and} rules}

{\bf Definition} A {\em strong multi-link} is a set of vertices $V$ and terminals
$T$ such that with play restricted to $V$, \short\ can
guarantee to cause $T$ to become a clique (in the reduced graph) if he plays first or second.
A {\em weak multi-link} is a set of vertices $V$ and terminals
$T$ such that with play restricted to $V$, \short\ can
guarantee to cause $T$ to become a clique if he plays first, but not if he plays
second. A {\em pivot} of a weak multi-link is any vertex which, if shorted,
makes the link become strong.

In either case the vertex set $V$ is said to be the {\em carrier} of the multi-link.
The definition of a strong multi-link includes the case where the carrier is
the empty set and the terminals are already a clique. Usually the terminals
$T$ are a subset of $\Gamma(V)$, but they do not need to be the whole of $\Gamma(V)$
so the existence of a strong multi-link $V,T$ does not necessarily imply that
$V$ is captured.

We can now state the generalized {\sc or}-rule:

\begin{theorem} {\em (`{\sc or}-rule')}:
If a set of weak multi-links between given terminals
has zero total intersection, then the union of the weak multi-links
is a strong multi-link between the terminals.
\end{theorem}

Proof: obvious, and can be compared with that of theorem \ref{th:winmulti}.

The `{\sc and}-rule' can be generalized in more than one way. The simplest
is, arguably:

\begin{theorem} {\em (`{\sc and}-rule')}
If two strong multi-links $V_1,T_1$ and $V_2,T_2$ have
intersecting terminals and non-intersecting carriers, then
the combination $V_1 \cup V_2, T_1 \oplus T_2$ is a weak multi-link,
pivoted by any vertex in $T_1 \cap T_2$.
\end{theorem}

Proof: use commutation and suppose the opening \short\ move is performed last.
The first strong multi-link guarantees that \short\ can cause
$T_1$ to become a clique, the second guarantees that \short\ can cause
$T_2$ to become a clique. When a vertex in $T_1 \cap T_2$ is then
shorted, $T_1 \oplus T_2$ becomes a clique.

As an example, we apply these rules to the multi-Shannon game shown in figure
\ref{fig:mgame_efgh}. $h$ is the carrier of a weak multi-link
from the left terminals to $f$; $g$ is the carrier of another such weak multi-link,
therefore by the {\sc or}-rule, $\{h,g\}$ is the carrier of a strong multi-link
between the left terminals and $f$. $f$ being adjacent to both right terminals,
there is a strong multi-link from $f$ to the right terminals with empty
carrier, and therefore by the {\sc and}-rule $\{h,g,f\}$ is the carrier of
a weak multi-link between all the terminals, pivoted by $f$. Finally, $e$ is the carrier of
another weak multi-link between all the terminals, so using the {\sc or}-rule a second
time, $\{h,g,f,e\}$ is the carrier of a strong multi-link between all the terminals
of the multi-Shannon game, so the game is won by \short.

\bibliographystyle{unsrt}
\bibliography{graphrefs}

\end{document}